\documentclass[12pt,twoside]{article}

\date{}

\begin{document}

\centerline{}

\centerline {\Large{\bf Fuzzy Anti-bounded Linear Operator}}

\centerline{}

\newcommand{\mvec}[1]{\mbox{\bfseries\itshape #1}}

\centerline{\bf {Bivas Dinda, T.K. Samanta and Iqbal H. Jebril}}

\centerline{}

\centerline{Department of Mathematics,}
\centerline{ Mahishamuri Ramkrishna
Vidyapith, West Bengal, India. }
\centerline{e-mail: bvsdinda@gmail.com}

\centerline{Department of Mathematics, Uluberia
College, West Bengal, India.}
\centerline{e-mail: mumpu$_{-}$tapas5@yahoo.co.in}

\centerline{Department of Mathematics, King Faisal University, Saudi Arabia. }
\centerline{e-mail: iqbal501@hotmail.com}

\centerline{}

\newtheorem{Theorem}{\quad Theorem}[section]

\newtheorem{definition}[Theorem]{\quad Definition}

\newtheorem{theorem}[Theorem]{\quad Theorem}

\newtheorem{remark}[Theorem]{\quad Remark}

\newtheorem{corollary}[Theorem]{\quad Corollary}

\newtheorem{note}[Theorem]{\quad Note}

\newtheorem{lemma}[Theorem]{\quad Lemma}

\newtheorem{example}[Theorem]{\quad Example}

\begin{abstract}
\textbf{\emph{Various types of fuzzy anti-continuity and fuzzy
anti-boundedness are defined. A few properties of them are
established. The intra and inter relation among various types of
fuzzy anti-continuity and fuzzy anti-boundedness are studied.}}
\end{abstract}


\section{Introduction}
Fuzzy set theory was first introduce by Zadeh\cite{zadeh} in 1965
and thereafter several authors contributed different articles on
this concept and applied on different branches of pure and applied
mathematics. The concept of fuzzy norm was introduced by Katsaras
\cite{Katsaras} in 1984. In 1992, Felbin\cite{Felbin1} introduced
the idea of fuzzy norm on a linear space. Cheng-Moderson
\cite{Shih-chuan} introduced another idea of fuzzy norm on a linear
space whose associated metric is same as the associated metric of
Kramosil-Michalek \cite{Kramosil}. Latter on Bag and Samanta
\cite{Bag1} modified the definition of fuzzy norm of
Cheng-Moderson \cite{Shih-chuan} and established the concept of
continuity and boundednes of a function with respect to
their fuzzy norm in \cite{Bag2}. \\\\
In this paper, various types of fuzzy anti-continuities and fuzzy
anti-boundedness; namely, fuzzy anti-continuity, sequential fuzzy
anti-continuity, strong fuzzy anti-continuity, weak fuzzy
anti-continuity, strong fuzzy anti-boundedness and weak fuzzy anti-boundedness are defined. The intra relations among fuzzy anti-continuities and intra relation
 among strongly fuzzy anti-bounded and weakly fuzzy anti-bounded are
 studied. Also, it is established a very simple criterion for fuzzy anti-continuity; namely, any linear operator between fuzzy anti-normed linear spaces is strongly and
 weakly fuzzy  anti-continuous if and only if it is strongly and
 weakly fuzzy  anti-bounded respectively.
\section{Preliminaries}
This section contain some basic definition and preliminary
results which will be needed in the sequel.\\
\begin{definition}
\cite{Schweizer}. A binary operation \, $\diamond \; : \; [\,0 \; ,
\; 1\,] \; \times \; [\,0 \; , \; 1\,] \;\, \longrightarrow \;\,
[\,0 \; , \; 1\,]$ \, is a\, $t$-conorm if
\,$\diamond$\, satisfies the
following conditions \, $:$ \\
$(\,i\,)\;\;$ \hspace{0.1cm} $\diamond$ \, is commutative and
associative ,
\\$(\,ii\,)$ \hspace{0.1cm} $a \;\diamond\;0 \;\,=\;\, a
\hspace{1.2cm}
\forall \;\; a \;\; \in\;\; [\,0 \;,\; 1\,]$ , \\
$(\,iii\,)$ \hspace{0.1cm} $a \;\diamond\; b \;\, \leq \;\, c
\;\diamond\; d$ \, whenever \, $a \;\leq\; c$  ,  $b \;\leq\; d$
 and  $a \, , \, b \, , \, c \, , \, d \;\; \in\;\;[\,0
\;,\; 1\,]$.
\end{definition}

\begin{remark}
\cite{Vijayabalaji}. $(\,a\,)\;$  For any  $\,r_{\,1} \; , \; r_{\,2}
\;\; \in\;\; (\,0 \;,\; 1\,)$ \, with \, $r_{\,1} \;>\;
r_{\,2}$ , there exist $r_{\,3} \;\in
\;(\,0 \;,\; 1\,)$ \, such that \, $r_{\,1} \;>\; r_{\,4} \;\diamond\;
r_{\,2}$ .
\\\\ $(\,b\,)\;$  For any  $\,r_{\,4} \;\,
\in\;\, (\,0 \;,\; 1\,)$ , there exist  $\,r_{\,5} \;\in\;\,(\,0 \;,\; 1\,)$ \,
such that \,and\, $r_{\,5}
\;\diamond\; r_{\,5} \;\leq\; r_{\,4}.$
\end{remark}

\begin{definition}\cite{dinda1}
Let $V$ be linear space over the field $F\,(=\,R\;or\;C\,)$.
A fuzzy subset $\;\nu\;$ of $\;V\,\times\,R\;$ is called a
\textbf{fuzzy antinorm} on $V$ if and only if for all
$x\,,\,y\;\in\;V$ and $c\;\in\;F$\\\\
$(i)$ For all $t\;\in\;R$ with $t\,\leq\,0\;,\;\nu(\,x\,,\,t\,)\;=\;1\;$;\\
$(ii)$ For all $t\;\in\;R$ with $t\,>\,0\;,\;\nu(\,x\,,\,t\,)\;=\;0\;$
if and only if $x\;=\;\theta$;\\
$(iii)$ For all $t\;\in\;R$ with $t\,>\,0\;,\;\nu\,(\,cx\,,\,t\,)=
\,\nu\,(\,x\,,\,\frac{t}{\mid c \mid}\,)\;$ if $c\neq\;0\,,c\;\in F$;\\
$(iv)$ For all $s\,,\,t\;\in\;R$ with $\;\nu\,(\,x\,+\,y\,,\,s\,+\,t\,)
\;\leq\;\nu\,(\,x\,,\,s\,)\;\diamond\;\nu\,(\,y\,,\,t\,)\;$;\\
$(v)$  $\mathop{\lim }\limits_{t\, \to\,\;\infty }\,\nu\,(\,x\,,\,t\,)=0.$
\end{definition}
We further assume that for any fuzzy anti-normed linear space
$\,(\,V\,,\,A^*\,),\\$
{\bf (vi)}$\;\;\,\nu\,(\,x\,,\,t\,)\,<\,1\;,\;\forall\,t>0\;
\Rightarrow\;x\,=\,\theta.\\ $
{\bf (vii)}$\;\,\nu(\,x\,,\,\cdot\,)\,$ is a continuous
function of $\mathbf{R}$ and strictly decreasing on the subset
 $\,\{\,t\,:\;0<\nu(x,t)<1\}\;$ of $\,\mathbf{R}.$

 \begin{theorem}\cite{dinda1}
Let $(\,V\,,\,A^*\,)$ be a fuzzy antinormed linear space
satisfying $(vi)$ and $(vii)$. Let $\left\| {\;x\;}
\right\|_{\,\alpha }^{\,\ast}\;=\;\wedge\,\{\,t\;:\;\nu\,(\,x\;,\;t\,)
\;\leq\;1\,-\,\alpha\}\;,\,\alpha\;\in\;(\,0\;,\;1\,).$ Also,
 let $\;\nu\,^\prime\,:\,V\,\times\,R\;\longrightarrow\;[\,0\,,\,1\,]$
 be defined by \[\nu\;^\prime\,(\,x\,,\,t\,)\;=\;\wedge\,\{\,1\,-\,
 \alpha\;:\;\left\| {\;x\;} \right\|_{\,\alpha }^{\,\ast}\;\leq\;t\}
 \;,\;if\;(\,x\,,\,t\,)\;\neq\;(\,\theta\,,\,0\,)\]
\[=\;1\hspace{4.4 cm},\;if\;(\,x\,,\,t\,)\;=\;(\,\theta\,,\,0\,)\hspace{-2 cm}\]
Then $\nu\,^\prime\;=\;\nu.$
\end{theorem}

\begin{definition}
\cite{Samanta1}. Let $(\,U\,,\,N^*\,)$ be a fuzzy antinormed
linear space. A sequence $\{x_n\}_n$ in $U$ is said to be
\textbf{convergent} to $x\,\in\,U$ if given $t\,>\,0\;,\;r\,
\in\,(\,0\,,\,1\,)$ there exist an integer $n_0\,\in\,\mathbf{N}$ such that
\[N^*\,(\;x_n\,-\,x\,,\,t\,)\;<\;r\;\;\;\forall\;n\,\geq\,n_0.\]
\end{definition}

\begin{definition}
\cite{Samanta1}. Let $(\,U\,,\,N^*\,)$ be a fuzzy antinormed
linear space. A sequence $\{x_n\}_n$ in $U$ is said to be
\textbf{cauchy sequence} to $x,\in\,U$ if given $t\,>\,0\;,
\;r\,\in\,(\,0\,,\,1\,)$ there exist an integer $n_0\,\in\,N$ such that
\[N^*\,(\;x_{n+p}\,-\,x_n\,,\,t\,)\;<\;r\;\;\;\forall\;n\,
\geq\,n_0\;\,,\,p=1,2,3,\cdots.\]
\end{definition}

\begin{definition}
\cite{Samanta1}. A subset $A$ of a fuzzy antinormed linear
space $(\,U\,,\,N^*\,)$ is said to be \textbf{bounded} if and
only if there exist $t\,>\,0\;,\;r\,\in\,(\,0\,,\,1\,)$ such
that \[N^*\,(\,x\,,\,t\,)\;<\;r\;\;\;\forall\;x\,\in\,A\]
\end{definition}

\section{Fuzzy Anti-Continuity  }

Let $(\,U\;,\;A^*\,)$ and $(\,V\;,\;B^*\,)$ be any two fuzzy
anti-normed linear spaces over the same field $F$.

\begin{definition}
A mapping $T\,:\,(\,U\;,\;A^\ast\,)\longrightarrow(\,V\;,\;B^\ast\,)$
is said to be {\bf fuzzy anti-continuous} at $x_0\;\in\;U$, if
for any given $\epsilon\,>\,0\,$,$\;\alpha\;\in\;(\,0\;,\;1\,)$
there exist $\delta\,=\,\delta\,(\,\alpha\,,\,\epsilon\,)\;>\;0\;$ ,
$\;\beta\,=\,\beta\,(\,\alpha\,,\,\epsilon\,)\;\in\;(\,0\;,\;1\,)$
such that for all $x\;\in\;U$
\[\nu_U\,(\;x\,-\,x_0\;,\;\delta\;)\;<\beta\;\;\Rightarrow\;\;\nu_V\,
(\;T(x)\,-\,T(x_0)\;,\;\epsilon\;)\;<\alpha\;.\]
\end{definition}

\begin{definition}
A mapping $T\,:\,(\,U\;,\;A^\ast\,)\longrightarrow(\,V\;,\;B^\ast\,)$
is said to be {\bf sequentially fuzzy anti-continuous} at $x_0\;\in\;U$,
if for any sequence $\{x_n\}_n\;,\;x_n\;\in\;U\;,\;\forall\;n$ with
$x_n\;\longrightarrow\;x_0$ implies $T(x_n)\;\longrightarrow\;T(x_0)$
in $V$, that is for all $t\,>\,0$,
\[\mathop {\lim }\limits_{n\; \to \;\infty }\nu_U\,(\;x_n\,-\,x_0\;,
\;t\;)\;=\;0\;\;\;\Rightarrow\;\;\;\mathop {\lim }\limits_{n\; \to
\;\infty }\nu_V\,(\;T(x_n)\,-\,T(x_0)\;,\;t\;)\;=\;0\;.\]
\end{definition}

\begin{definition}
A mapping $T\,:\,(\,U\;,\;A^\ast\,)\longrightarrow(\,V\;,\;B^\ast\,)$
is said to be {\bf strongly fuzzy anti-continuous} at $x_0\;\in\;U$,
if for any given $\epsilon\,>\,0\,$ there exist $\delta\,=\,\delta\,
(\,\alpha\,,\,\epsilon\,)\;>\;0\;$ such that for all $x\;\in\;U$,
\[\nu_V\,(\;T(x)\,-\,T(x_0)\;,\;\epsilon\;)\;\leq\;\nu_U\,(\;x\,-\,x_0\;,\;\delta\;)\]
\end{definition}

\begin{definition}
A mapping $T\,:\,(\,U\;,\;A^\ast\,)\longrightarrow(\,V\;,\;B^\ast\,)$
is said to be {\bf weakly fuzzy anti-continuous} at $x_0\;\in\;U$,
if for any given $\epsilon\,>\,0\;,\;\alpha\;\in\;(\,0\,,\,1\,)$
there exist $\delta\,=\,\delta\,(\,\alpha\,,\,\epsilon\,)\;>\;0\;$
such that for all $x\;\in\;U$,
\[\nu_U\,(\;x\,-\,x_0\;,\;\delta\;)\;\leq\;1\,-\,\alpha\;\;
\Rightarrow\;\;\nu_V\,(\;T(x)\,-\,T(x_0)\;,\;\epsilon\;)\;
\leq\;1\,-\,\alpha\;.\]
\end{definition}

\begin{theorem}
If a mapping $T$ from a fuzzy anti-normed linear space
$(\,U\,,\,A^*\,)$ to a fuzzy anti-normed linear space $(\,V\,,\,B^*\,)$
is strongly fuzzy anti-continuous then it is weakly fuzzy anti-continuous.
But not conversely.
\end{theorem}
{\bf Proof.}$\;\;$
Obvious.

To show the converse result may not be true we consider the following example.
\begin{example}
As in the example of Note 3.3 of \cite{dinda}, we consider the fuzzy
anti-normed linear spaces $(\,X\,,\,\nu_1\,)$ and $(\,X\,,\,\nu_2\,)$.
Let $f(x)\,=\,\frac{x^4}{1+x^4}\;\;\forall\,x\,\in\,\mathbf{R}.$ Now
from Example 3 of \cite{Samanta} it directly follows that $f$ is not
strongly fuzzy anti-continuous. Here we now show that $f$ is weakly
fuzzy anti-continuous on $X$.\\
Let $x_0\,\in\,X\;,\;\epsilon>0$ and $\delta\,\in\,(0,1)$. Now\\
$\nu_2(\,f(x)-f(x_0)\,,\,\epsilon\,)\;<\;1-\alpha\;\;$ if $\;
\;\frac{k\mid f(x)-f(x_0\mid}{\epsilon + k\mid f(x)-f(x_0\mid}\;
<\;1-\alpha\\$
i.e., if \[\frac{\epsilon}{\epsilon \,+ k\,\mid \,\frac{x^4}{1+x^2}
\;-\;\frac{x_0^4}{1\,+\,x_0^2}\, \mid}\;\geq\;\alpha\]
i.e., if \[\frac{\frac{\epsilon\,(1\,+\,x^2)\,(1\,+\,x_0^2)}{k\,\mid\,
x\,+\,x_0\,\mid\;\mid \,x^2\, x_0^2 \,+\,x^2\,+\, x_0^2\,
\mid}}{\frac{\epsilon\,(1\,+\,x^2)\,(1\,+\,x_0^2)}{k\,\mid\,x\,
+\,x_0\,\mid\;\mid \,x^2\,x_0^2\,+\,x^2\,+\,x_0^2\,\mid}\,+\,
\mid x-x_0\,\mid}\;\geq\;\alpha\]
i.e., if \[\alpha\,\mid\,x\,-\,x_0\,\mid\;\leq\;(1-\alpha)\,
\frac{\epsilon}{k}\;\frac{(1\,+\,x^2)\,(1\,+\,x_0^2)}{\mid\,x\,
+\,x_0\,\mid\;\mid \,x^2\,x_0^2\,+\,x^2\,+\,x_0^2\,\mid}\]
\[\leq\;(1-\alpha)\,\frac{\epsilon}{k}\;\hspace{3 cm}\]
So, depending upon $\;(1-\alpha)\,\frac{\epsilon}{k}\;$ we may
choose $\;\delta\,>\,0\;$ such that \\
$\alpha\,(\,\delta\,+\,\mid\,x\,-\,x_0\,\mid)\;\leq\;\delta\;\;
\;\;i.e.,\;\;\nu_1(\,x\,-\,x_0\,,\,\delta\,)\;<\;1-\alpha.\\$
Thus we see that for every $\epsilon\,>\,0\,,\,\alpha\,\in\,(0,1)\;
\;\exists\;\delta>0\;$ such that
\[\nu_1(\,x\,-\,x_0\,,\,\delta\,)\;<\;1-\alpha\;\Rightarrow\;
\nu_2(\,f(x)\,-\,f(x_0)\,,\,\epsilon\,)\;<\;1-\alpha.\]
i.e., $f$ is weakly fuzzy anti-continuous at $x_0$.
\end{example}

\begin{theorem}\label{th52}
 A mapping $T$ from a fuzzy anti-normed linear space $(\,U\,,\,A^*\,)$
 to a fuzzy anti-normed linear space $(\,V\,,\,B^*\,)$ is fuzzy
 anti-continuous if and only if it is sequentially fuzzy anti-continuous.
\end{theorem}
{\bf Proof.}$\;\;$
The proof of the above theorem is directly follows from Theorem 13 of \cite{Samanta}.

\begin{theorem}
If a mapping $T$ from a fuzzy anti-normed linear space $(\,U\,,\,A^*\,)$
to a fuzzy anti-normed linear space $(\,V\,,\,B^*\,)$ is strongly fuzzy
anti-continuous then it is sequentially fuzzy anti-continuous.
\end{theorem}
{\bf Proof.}$\;\;$
The proof of the above theorem is directly follows from Theorem 12 of \cite{Samanta}.

\begin{theorem}
Let $T\,:\,(\,U\;,\;A^\ast\,)\longrightarrow(\,V\;,\;B^\ast\,)$ be
a linear operator. If $T$ is sequentially fuzzy anti-continuous at
a point $x_0\;\in\;U\,$, then it is sequentially fuzzy anti-continuous on $U\,.$
\end{theorem}
{\bf Proof.$\;\;$}
Let $x\;\in\;U$ be an arbitrary point and let $\{x_n\}_n$ be a
sequence in $U$ such that $x_n\;\longrightarrow\;x$. Then $\forall\;t\;>\;0$\\
\[\mathop {\lim }\limits_{n\;\,\to \,\;\infty }\,\nu\,_U\,(\;x_n\,
-\,x\,,\,t\,)\;=\;0\hspace{7 cm}\]
\[i.e.,\;\mathop {\lim }\limits_{n\;\,\to \,\;\infty }\,\nu\,_U\,(\;(\,x_n\,
-\,x\,+\,x_0\,)\,-\,x_0\,,\,t\,)\;=\;0\hspace{6 cm}\]
Since $T$ is sequentially fuzzy anti-continuous at
$x_0\;\;\forall\;t\;>\;0\;$ we have \\
\[\mathop {\lim }\limits_{n\;\,\to \,\;\infty }\,\nu\,_V\,(\;
(\,x_n\,-\,x\,+\,x_0\,)\,-\,x_0\,,\,t\,)\;=\;0\hspace{4.5 cm}\]
\[i.e.,\;\mathop {\lim }\limits_{n\;\,\to\,\;\infty}\,\nu\,_V\,
(\;T(\,x_n\,)\,-\,T(\,x\,)\,+\,T(\,x_0\,)\,-\,T(\,x_0\,)\;,\,t\;)
\;=\;0\hspace{6 cm}\]
\[i.e.,\;\mathop {\lim }\limits_{n\;\,\to \,\;\infty }\,\nu\,
_V\,(\;T(\,x_n\,)\,-\,T(\,x\,)\;,\,t\;)\;=\;0\hspace{7 cm}\]
Thus, \[\mathop {\lim }\limits_{n\;\,\to \,\;\infty }\,\nu\,_U\,
(\;x_n\,-\,x\,,\,t\,)\;=\;0\;,\;\forall\;t\;>\;0\;\Rightarrow\;
\mathop {\lim }\limits_{n\;\,\to \,\;\infty }\,\nu\,_V\,(\;T(\,x_n\,)
\,-\,T(\,x\,)\;,\,t\;)\;=\;0\;,\;\forall\;t\;>\;0.\]
Hence the proof.
\bigskip

\section{Fuzzy Anti-Boundedness }

\begin{definition}
A mapping $T\,:\,(\,U\;,\;A^\ast\,)\longrightarrow(\,V\;,\;B^\ast\,)$
is said to be \textbf{strongly fuzzy anti-bounded} on $U$ if and only
if there exist a positive real number $M$ such that for all $x\,\in\,U$
and for all $t\,\in\,R^+$,
\[\nu_V\,(\;T(x)\,,\,t\,)\;\leq\;\nu_U\,(\,x\,,\,\frac{t}{M}\,)\]
\end{definition}

\begin{example}
The zero and identity operators are strongly fuzzy anti-bounded.
\end{example}

\begin{example}
It is an example of a strongly fuzzy anti-bounded linear operator
other than the zero and identity operator. \\
Let $(\,V\,,\,\|.\|\,)$ be a normed linear space over the field
 $\;K\,(=\,R\;or\;C).$  Let $\alpha_1,\alpha_2\,\in\,R$ such that
  $\alpha_1>\alpha_2>0$. Again, let $\nu_1,\nu_2\,:\,V\times R^+
  \longrightarrow\,[0,1]$ be defined by\\
\[\nu_1\,(x,t)=\frac{\alpha_1\|x\|}{t+\alpha_1\|x\|}\;\;and \;\;
\nu_2\,(x,t)=\frac{\alpha_2\|x\|}{t+\alpha_2\|x\|}\]
Also, define  $a\,\diamond\,b,\,=\,max\{a,b\}$ for all $a,b\,\in\,[0,1].$\\
Now we shall first show that $(V,\nu_1)$ and $(V,\nu_2)$ are fuzzy
anti-normed linear space.\\
$(i)\;\;$ The condition (i) is obvious.\\
$(ii)\;\;$ $\nu_1(x,t)=0\;\Leftrightarrow\;
\frac{\alpha_1\|x\|}{t+\alpha_1\|x\|}=0\;\Leftrightarrow\;
\|x\|=0\;\Leftrightarrow\; x=\theta.$\\\\
$(iii)\;\;$ Let $c\,\in\,K\,$ and $c\,\neq\,0$
\[\nu_1(\,cx\,,\,t\,)\;=\;\frac{\alpha_1\,\parallel\,cx\,
\parallel}{t\,+\,\alpha_1\,\parallel\,cx\,\parallel}\]
\[=\;\frac{\alpha_1\,\parallel\,x\,\parallel}{\frac{t}{\mid\,c\,
\mid}\,+\,\alpha_1\,\parallel\,x\,\parallel}\;=\;\nu_1(\,x\,,\,
\frac{t}{\mid\,c\,\mid})
\hspace{-5 cm}\]
$(iv)\;\;$ \[\nu_1(\,x+y\,,\,s+t\,)\;=\;\frac{\alpha_1\,\parallel\,
x+y\,\parallel}{s\,+\,t\,+\,\alpha_1\,\parallel\,x+y\,\parallel}\]
\[=\;\frac{1}{\frac{s\,+\,t}{\alpha_1\,\parallel\,x+y\,
\parallel}\;+\;1}{\hspace{-1.5 cm}}\]
\[\leq\;\frac{1}{\frac{s\,+\,t}{\alpha_1\,\parallel\,x\,\parallel\,
+\,\alpha_1\,\parallel\,y\,\parallel}\;+\;1}{\hspace{-2.5 cm}}\]
\[=\;\frac{\alpha_1\,\parallel\,x\,\parallel\,+\,\alpha_1\,
\parallel\,y\,\parallel}{s\,+\,t\,+\,\alpha_1\,\parallel\,x\,\parallel\,
+\,\alpha_1\,\parallel\,y\,\parallel}{\hspace{-4.5 cm}}\]
Now if \[\nu_1(x,s)\;\geq\;\nu_1(y,t)\;\Rightarrow
\;\frac{\alpha_1\,\parallel\,x\,\parallel}{s\,+\,\alpha_1\,
\parallel\,x\,\parallel} \;\geq\frac{\alpha_1\,\parallel\,y\,
\parallel}{t\,+\,\alpha_1\,\parallel\,y\,\parallel};\]
\[\;\Rightarrow\;t\,\parallel\,x\,\parallel\,-\,s\,\parallel\,y\,
\parallel{\hspace{-1 cm}}\]
Therefore, \[\;\;\frac{\alpha_1\,\parallel\,x\,\parallel\,+\,
\alpha_1\,\parallel\,y\,\parallel}{s\,+\,t\,+\,\alpha_1\,
\parallel\,x\,\parallel\,
+\,\alpha_1\,\parallel\,y\,\parallel}\;-\;\frac{\alpha_1\,
\parallel\,x\,\parallel}{s\,+\,\alpha_1\,\parallel\,x\,\parallel}\;\leq\;0.\]

Thus \[\nu_1(x+y\,,\,s+t)\;\leq\;\frac{\alpha_1\,\parallel\,x\,
\parallel\,+\,\alpha_1\,\parallel\,y\,\parallel}{s\,+\,t\,
+\,\alpha_1\,\parallel\,x\,
\parallel\,+\,\alpha_1\,\parallel\,y\,\parallel}\]
\[\;\leq\;\frac{\alpha_1\,\parallel\,x\,\parallel}{s\,+
\,\alpha_1\,\parallel\,x\,\parallel}\; =\;\nu_1(\,x\,,\,s\,)
\diamond\nu_1(\,y\,,\,t\,){\hspace{-4 cm}}\]
Again if $\nu_1(y,t)\;\geq\;\nu_1(x,s)$ Similarly it can be shown that
\[\nu_1(x+y\,,\,s+t)\;\leq\;\frac{\alpha_1\,\parallel
\,y\parallel}{t\,+\,\alpha_1\,\parallel\,y\parallel}\;
=\;\nu_1(\,x\,,\,s\,)\diamond\nu_1(\,y\,,\,t\,)\]
Hence  \[\nu_1(x+y\,,\,s+t)\;\leq\;\nu_1(\,x\,,\,s\,)
\diamond\nu_1(\,y\,,\,t\,)\]
$(v)\;\;$\[\mathop {\lim }\limits_{t\; \to \;\infty }\,
\nu_1(x,t)\;=\;\mathop {\lim }\limits_{t\; \to \;\infty }
\frac{\alpha_1\,\parallel\,x\,\parallel}{t\,+\,\alpha_1\,
\parallel\,x\,\parallel}\;=\;0\]
Hence $(\,V\,,\,\nu_1\,)$ is a fuzzy anti-normed linear space.
Similarly $(\,V\,,\,\nu_2\,)$ is also fuzzy anti-normed linear space.\\\\

We now define a mapping $\;T\;:\;(\,V\,,\,\nu_1\,)
\;\longrightarrow\;(\,V\,,\,\nu_2\,)\;$ by $\;T(x)\,=\,rx\;$
 where $r\,(\neq\,0)\,\in\,\mathbf{R}$ is fixed. Clearly $T$
 is a linear operator.\\
Let us choose an arbitrary but fixed $\;M\,>\,0\;$
such that $\,M\;\geq\;\mid\,r\,\mid\,$ and $\;x\,\in\,V.\;$
Now,\[\,M\;\geq\;\mid\,r\,\mid\;\Rightarrow\;\alpha_1\,M\,
\parallel\,x\,\parallel\;\geq\;\alpha_2\,\mid\,r\,\mid\,
\parallel\,x\,\parallel\hspace{6 cm}\]
\[\Rightarrow\;t\,+\,\alpha_1\,M\,\parallel\,x\,\parallel\;
\geq\;t\,+\,\alpha_2\,\mid\,r\,\mid\,\parallel\,x\,\parallel\;
\;\;\;\;\forall\;t>0.\]
\[\Rightarrow\;\frac{t}{t\,+\,\alpha_2\,\mid\,r\,\mid\,\parallel\,
x\,\parallel}\;\geq\;\frac{t}{t\,+\,\alpha_1\,M\,\parallel\,x\,\parallel}\;\;
\;\;\;\forall\;t>0.\]
\[\Rightarrow\;\frac{t}{t\,+\,\alpha_2\,\parallel\,rx\,\parallel}
\;\geq\;\frac{\frac{t}{M}}{\frac{t}{M}\,+\,\alpha_1\,\parallel\,x\,
\parallel}\;\;\;\;\;\forall\;t>0.\hspace{1 cm}\]
\[\Rightarrow\;1\,-\,\frac{t}{t\,+\,\alpha_2\,\parallel\,rx\,
\parallel}\;\leq\;1\,-\,\frac{\frac{t}{M}}{\frac{t}{M}\,+\,
\alpha_1\,\parallel\,x\,
\parallel}\;\;\;\;\;\forall\;t>0.\hspace{-.75 cm}\]
\[\Rightarrow\;\frac{\alpha_2\,\parallel\,rx\,\parallel}{t\,+\,
\alpha_2\,\parallel\,rx\,\parallel}\leq\;\frac{\alpha_1\,
\parallel\,x\,\parallel}
{\frac{t}{M}\,+\,\alpha_1\,\parallel\,x\,\parallel}\;\;\;\;
\;\forall\;t>0.\hspace{1.25 cm}\]
i.e.,\[\nu_2(\,T(x)\,,\,t\,)\;\leq\;\nu_1(\,x\,,\,\frac{t}{M})\;
\;\;\;\;\forall\;t>0\;\;and\;\;\forall\;x\,\in\,V.\]
Hence $T$ is strongly fuzzy anti-bounded on $V$.
\end{example}

\begin{definition}
A mapping $T\,:\,(\,U\;,\;A^\ast\,)\longrightarrow
(\,V\;,\;B^\ast\,)$ is said to be \textbf{weakly fuzzy anti-bounded}
 on $U$ if and only if for any $\alpha\,\in\,(\,0\,,\,1\,)$
 there exist $M_\alpha\,(\,>\,0)$ such that for all $x\,\in\,U$
 and for all $t\,\in\,R^+$,
\[\nu_U\,(\,x\,,\,\frac{t}{M_\alpha}\,)\;\leq\;1\,-\,\alpha\;
\;\Rightarrow\;\nu_V\,(\;T(x)\,,\,t\,)\; \leq\;1\,-\,\alpha\;\;\]
\end{definition}

\begin{theorem}
Let $T\,:\,(\,U\;,\;A^\ast\,)\longrightarrow(\,V\;,\;B^\ast\,)$
 be a linear operator. If $T$ is strongly fuzzy anti-bounded then
 it is weakly fuzzy anti-bounded. But not conversely.
\end{theorem}
{\bf Proof.}$\;\;$
First we suppose that $T$ is strongly fuzzy anti-bounded.
Then there exist $\,M\,>\,0\,$ such that $\;\forall\,x\,\in\,U\;$
and $\;\forall\,t\,\in\,R\;$, \[\nu_V\,(\;T(x)\,,\,t\,)
\;\leq\;\nu_U\,(\,x\,,\,\frac{t}{M}\,)\]
Thus for any $\alpha\,\in\,(\,0\,,\,1\,)$, there exists
$M_\alpha\;(\,=\,M\,)$ such that \[\nu_U\,(\,x\,,
\,\frac{t}{M_\alpha}\,)\;\leq\;1\,-\,\alpha\;\;\Rightarrow
\;\nu_V\,(\;T(x)\,,\,t\,)\; \leq\;1\,-\,\alpha\;\;\]
Hence $T$ is weakly fuzzy anti-bounded.

The converse of the above theorem is not necessarily true. For example
\begin{example}
Let $(\,V\,,\,\|\cdot\|\,)\;$ be a linear space over the
field $\;K(=R\;or\;C)$ and $\;\nu_1,\;\nu_2\;:\;V\times R
\;\longrightarrow\;[0,1]$ be defined by
\[\nu_1(x,t)=\frac{2 \|x\|^2}{t^2+\|x\|^2}\;\;\;,\;if\;t>\|x\|\]
\[ =\,1\hspace{2 cm},\;if\;0<t\leq\|x\| \hspace{-2.2 cm}\]
\[and \hspace{3 cm} \nu_2(x,t)=\frac{\|x\|}{t+\|x\|}\hspace{6.7 cm}\]
Also define $\;a\diamond b=\;max\{a,b\}\\$
Already we have seen that $\;(\,V\,,\,\nu_2\,)\;$ is a fuzzy
anti-normed linear space. Now we shall prove that
$\;(\,V\,,\,\nu_1\,)\;$ is a fuzzy anti-normed linear space.\\\\
$(i)\;\;$ The condition (i) is obvious.\\
$(ii)\;\;$ $\nu_1(x,t)=0\;\Leftrightarrow\;
\frac{2 \|x\|^2}{t^2+\|x\|^2}=0\;\Leftrightarrow\;\|x\|=0
\;\Leftrightarrow\;x=\theta.\\$
$(iii)\;\;$ Let, $\;c\,\in K\;$ and $\;c\neq 0.\;$ If $\;t>\|cx\|\,,$
\[\nu_1(cx,t)\;=\;\frac{2 \|cx\|^2}{t^2+\|cx\|^2}\;=
\;\frac{2 |c|^2\|x\|^2}{t^2+|c|^2\|x\|^2}\;=\;
\frac{2 \|x\|^2}{(\frac{t}{|c|})^2+\|x\|^2}\;=\;\nu_1(x,\frac{t}{|\,c|})\]
Again if $\;0<t\leq\|cx\|$ then $\nu_1(cx,t)=1\\$
and $\;0<t\leq\|cx\|\;\Rightarrow\;0<\frac{t}{|c|}\leq \|x\|
\;\Rightarrow\;\nu_1(x,\frac{t}{|c|})=1$\\
$(iv)\;\;$Let $s,t\,\in\,\mathbf{R}^+\;,\;x,y\,\in\,V\\$
If $\;\;0\;<\;s+t\;\leq\;\|x+y\|$, we have the following possibilities
\[(a)\;0\;<\;s\;\leq\;\|x\|\;\;and\;\;0\;<\;t\;\leq\;\|y\|\hspace{-1 cm}\]
\[(b)\;0\;<\;s\;\leq\;\|x\|\;\;and\;\;t\;>\|y\|\]
\[(c)\;0\;<\;t\;\leq\;\|y\|\;\;and\;\;s\;>\;\|x\|\]
In each case $\nu_1(\,x+y\,,\,s+t\,)\;=\;1\;=
\nu_1(\,x\,,\,s\,)\diamond \nu_1(\,y\,,\,t\,)$
Again, if $\;s+t\;>\;\|x+y\|$, we have the following four possibilities
\[(a)\;s\;>\;\|x\|\;,\;t\;\leq\;\|y\|\]
\[(b)\;s\;\leq\;\|x\|\;,\;t\;>\;\|y\|\]
\[(c)\;s\;\leq\;\|x\|\;,\;t\;\leq\;\|y\|\]
\[(d)\;s\;>\;\|x\|\;,\;t\;>\;\|y\|\]
In the cases (a), (b), (c)
\[\nu_1(\,x+y\,,\,s+t\,)\;=\;\frac{2\,\|x+y\|^2}{(s+t)^2\,+\,\|x+y\|^2}\]
\[<\;1\;=\;\nu_1(\,x\,,\,s\,)\diamond \nu_1(\,y\,,\,t\,)\hspace{-3.75 cm}\]
So, we now suppose that $\;s\;>\;\|x\|\;\;and \;\;t\;>\;\|y\|$
Now, $\;s+t\;>\;\|x\|\;+\;\|y\|\;\geq\;\|x+y\|$. Therefore,
\[\nu_1(\,x+y\,,\,s+t\,)\;=\;\frac{2\,\|x+y\|^2}{(s+t)^2\,+\,\|x+y\|^2}\]
\[\leq\;\frac{2(\;\|x\|\;+\;\|y\|\;)^2}{(s+t)^2\;+\;(\;\|x\|\;+\;\|y\|\;)^2}
\hspace{-4.75 cm}\]
Hence we have
\[\nu_1(\,x+y\,,\,s+t\,)\;\leq\;\frac{2(\;\|x\|\;+
\;\|y\|\;)^2}{(s+t)^2\;+\;(\;\|x\|\;+\;\|y\|\;)^2}\]
\[\leq\;\frac{2\,\|y\|^2}{t^2\;+\;\|y\|^2}\;=\;\nu_1(\,y\,,\,t\,)\hspace{-2.75 cm}\]
when $\;\;\nu_1(x,s)\;\leq\;\nu_1(y,t)\\$
Similarly,
\[\nu_1(\,x+y\,,\,s+t\,)\;\leq\;\frac{2(\;\|x\|\;+\;\|y\|
\;)^2}{(s+t)^2\;+\;(\;\|x\|\;+\;\|y\|\;)^2}\]
\[\leq\;\frac{2\,\|x\|^2}{t^2\;+\;\|x\|^2}\;=\;\nu_1(\,x\,,\,s\,)\hspace{-2.75 cm}\]
when $\;\;\nu_1(y,t)\;\leq\;\nu_1(x,s)\\$
Thus \[\nu_1(\,x+y\,,\,s+t\,)\;\leq\;\nu_1(\,x\,,\,s\,)\diamond \nu_1(\,y\,,\,t\,)\]

\[(v)\;\;\mathop {\lim }\limits_{t\; \to \;\infty }\,
\nu_1(x,t)\;=\;\mathop {\lim }\limits_{t\; \to \;\infty }\,
\frac{2 \|x\|^2}{t^2+\|x\|^2}=0 \hspace{6.4 cm}\]
Thus we see that $\;(\,V\,,\,\nu_1\,)$ is a fuzzy anti-normad linear space.\\\\
Now we define a linear operator $\;T\,:\,(\,U\,,\,\nu_1\,)
\longrightarrow (\,V\,,\,\nu_2\,)\;$ by $\;T(x)=x\;\;,\;
\forall\;x\,\in \,V.\;$ Let, $\alpha\,\in\,(0,1)\;,\; x\,
\in\,V\;and t\,\in\,R^+\;$ and choose $\;M_\alpha=
\frac{1}{1-\alpha}.$ We now prove that
\[\nu_1(\,x\,,\,\frac{t}{M_\alpha})\;\leq\;1-\alpha\;\Rightarrow
\;\nu_2(\,T(x)\,,\,t\,)\;\leq\;1-\alpha.\]
\[\nu_1(\,x\,,\,\frac{t}{M_\alpha})\;\leq\;1-\alpha \hspace{10 cm}\]
\[\Rightarrow\;\frac{2\,{\parallel\,x\,\parallel}^2}{t^2\,(1-\alpha)^2
\;+\;{\parallel\,x\,\parallel}^2}\;\leq\;1-\alpha\hspace{10 cm}\]
 \[\Rightarrow\;1\;-\;\frac{2\,{\parallel\,x\,\parallel}^2}{t^2\,
 (1-\alpha)^2\;+\;{\parallel\,x\,\parallel}^2}\;\geq\;1-(1-\alpha)\;
 =\;\alpha\hspace{10 cm}\]
\[\Rightarrow\;\frac{t^2\,(1-\alpha)^2\;-\;{\parallel\,x\,
\parallel}^2}{t^2\,(1-\alpha)^2\;+\;{\parallel\,x\,\parallel}^2}
\;\geq\;\alpha \hspace{10 cm}\] \[\Rightarrow\;t^2\,(1-\alpha)^3
\;\geq\;(1+\alpha)\,{\parallel\,x\,\parallel}^2  \hspace{10 cm}\]
\[\Rightarrow\;\parallel\,x\,\parallel\;\leq\;\frac{t(1-\alpha)\,
\sqrt{1-\alpha}}{\sqrt{1+\alpha}} \hspace{10 cm}\] \[\Rightarrow
\;t\,+\,\parallel\,x\,\parallel\;\leq\;t\,\frac{(1-\alpha)\,
\sqrt{1-\alpha}\,+\,\sqrt{1+\alpha}}{\sqrt{1+\alpha}} \hspace{10 cm}\]
\[\Rightarrow\;\frac{t}{t\,+\,\parallel\,x\,\parallel}\;\geq
\;\frac{\sqrt{1+\alpha}}{(1-\alpha)\,\sqrt{1-\alpha}\,+
\,\sqrt{1+\alpha}} \hspace{10 cm}\] \[\Rightarrow\;1\;-
\;\frac{t}{t\,+\,\parallel\,x\,\parallel}\;\leq\;1\;-\;
\frac{\sqrt{1+\alpha}}{(1-\alpha)\,\sqrt{1-\alpha}\,+\,
\sqrt{1+\alpha}}
\hspace{10 cm}\] \[\Rightarrow\;\frac{\parallel\,x\,\parallel}{t\,+\,
\parallel\,x\,\parallel}\;\leq\;\frac{(1-\alpha)\sqrt{1-\alpha}}{(1-
\alpha)\sqrt{1-\alpha} \,+\,\sqrt{1+\alpha}}\hspace{10 cm}\]
 Now
 \[\frac{(1-\alpha)\sqrt{1-\alpha}}{(1-\alpha)\sqrt{1-\alpha} \,+
 \,\sqrt{1+\alpha}}\;\leq\;1-\alpha\Leftrightarrow\;\sqrt{1-\alpha}
 \;\leq\;(1-\alpha)\sqrt{1-\alpha} \,+\,\sqrt{1+\alpha}\]
 \[\Leftrightarrow\;\alpha\,\sqrt{1-\alpha}\leq\,\sqrt{1+\alpha}\hspace{-3.5 cm}\]
 \[\Leftrightarrow\;1\,+\,\alpha\,+\,\alpha^3\;\geq\;\alpha^2\hspace{-3.25 cm}\]
 which is true for all $\alpha\,\in\,(0,1).\\$
 Hence  \[\nu_1(\,x\,,\,\frac{t}{M_\alpha})\;\leq\;1-\alpha
 \;\Rightarrow\;\nu_2(\,T(x)\,,\,t\,)\;\leq\;1-\alpha.\]
 Thus  $T$ is weakly fuzzy anti-bounded on $X$.
 \end{example}

\begin{definition}
A linear operator $T\,:\,(\,U\;,\;A^\ast\,)
\longrightarrow(\,V\;,\;B^\ast\,)$ is said to be
\textbf{uniformly fuzzy anti-bounded} if and only
if there exist $M\,>\,0$ such that \[\;\left\| {\;T(x)\;}
\right\|_{\,\alpha }^{\,\ast}\;\geq\;M\,\left\| {\;x\;}
\right\|_{\,\alpha}^{\,\ast}\;\;\;,\;\alpha\,\in\,(\,0\,,\,1\,)\]
where  $\{\;\left\| {\;\cdot\;} \right\|_{\,\alpha }^{\,\ast}\;
:\;\alpha\,\in\,(\,0\,,\,1\,)\}$ is ascending family of fuzzy $\alpha$-anti-norms.
\end{definition}

\begin{theorem}
Let $T\,:\,(\,U\;,\;A^\ast\,)\longrightarrow(\,V\;,\;B^\ast\,)$
be a linear operator and $(\,U\;,\;A^\ast\,)$ and $(\,V\;,\;B^\ast\,)$
satisfies (vi) and (vii). Then  $T$ is strongly fuzzy anti-bounded
if and only if it is uniformly fuzzy anti-bounded with respect to
fuzzy $\alpha$-anti-norms.
\end{theorem}
{\bf Proof.}$\;\;$
Let $\{\;\left\| {\;\cdot\;} \right\|_{\,\alpha }^{\,\ast}\;:
\;\alpha\;\in\;(\,0\;,\;1\,)\;\}$ be ascending family of
$\alpha$-norms. First suppose that $\,T\,$ is strongly fuzzy
anti-bounded. Then there exist $\,M\,>\,0\,$ such that
 $\;\forall\,x\,\in\,U\;$ and $\;\forall\,s\,\in\,R\;$,
 \[\nu_V\,(\;T(x)\,,\,t\,)\;\leq\;\nu_U\,(\,x\,,\,\frac{s}{M}\,)
 \;\;\;i.e.,\,\nu_V\,(\;T(x)\,,\,t\,)\;\leq\; \nu_U\,(\,M\,x\,,\,s\,)\]
$\left\| {\;M\,x\;} \right\|_{\,\alpha }^{\,\ast}\;>\;t\;\Rightarrow
\;\wedge\,\{\,s\;:\;\nu\,(\,M\,x\;,\;s\,)\;\leq\;1\,-\,\alpha\}\;>\;t\;.$
\[\;\Rightarrow\;\exists\;s_0\,>\,t\;such \;that\;\nu\,(\,M\,x\;,\;s_0\,)
\;\leq\;1\,-\,\alpha \hspace{1 cm}\]
\[\;\Rightarrow\;\exists\;s_0\,>\,t\;such \;that\;\nu\,(\,T\,(\,x\,)\;,
\;s_0\,)\;\leq\;1\,-\,\alpha \hspace{.6 cm}\]
\[\;\Rightarrow\;\left\| {\;T\,(x)\;} \right\|_{\,\alpha }^{\,\ast}
\;\geq\;s_0\;>\;t\;\hspace{4.2 cm}\]
Hence $\;\left\| {\;T\,(x)\;} \right\|_{\,\alpha }^{\,\ast}\;\geq
\;\left\| {\;M\,x\;} \right\|_{\,\alpha }^{\,\ast}\;=\;M\,\left\|
{\;x\;} \right\|_{\,\alpha }^{\,\ast}\;.$\\
Thus $T$ is uniformly fuzzy anti-bounded.
\\\\
Conversely, suppose that there exist $M\,>\,0$ such that
$\forall\;x\;\in\;U$ and $\forall\;\alpha\;\in\;(\,0\,,\,1\,)$
 \[\;\left\| {\;T(x)\;} \right\|_{\,\alpha }^{\,\ast}\;\geq\;M\,
 \left\| {\;x\;} \right\|_{\,\alpha}^{\,\ast}\;\;\]
Let, $\;p\;>\;\nu_U\,(\,Mx\,,\,s\,)\;\Rightarrow\;\;p\;>\;\wedge\,
\{\;\alpha\,\in\;(\,0\,,\,1\,)\,:\,\left\| {\;M\,x\;}\right\|_{\,
\alpha }^{\,\ast}\;\leq\;s\;\}\\ \Rightarrow\;\;$ There exist
$\;\alpha_0\,\in\;(\,0\,,\,1\,)\,$ such that $\,p\,>\,\alpha_0\;
\,and\;\left\| {\;M\,x\;} \right\|_{\,\alpha }^{\,\ast}\;\leq\;s\;\\
\Rightarrow\;\;\left\| {\;T\,(x)\;} \right\|_{\,\alpha }^{\,\ast}
\;\leq\;s\\ \Rightarrow\;\;\nu_V\,(\,T(x)\,,\,s\,)\;\leq\;1\,-\,
\alpha_0\;<\;p.\\ $ Hence, $\;\nu_V\,(\,T(x)\,,\,s\,)\;\leq\;
\nu_U\,(\,Mx\,,\,s\,)\;=\;\nu_U\,(\,x\,,\,\frac{s}{M}\,).$\\
Thus $T$ is strongly fuzzy anti-bounded.

\begin{theorem}
Let $T\,:\,(\,U\;,\;A^\ast\,)\longrightarrow(\,V\;,\;B^\ast\,)$
be a linear operator. Then,\\
{\bf(i)$\;$} $T$ is strongly fuzzy anti-continuous on $U$ if
$\;T\;$ is strongly fuzzy anti-continuous at a point $x_0\;\in\;U\,.$\\
{\bf(ii)$\;$} $T$ is strongly fuzzy anti-continuous if and only if
$T$ is strongly fuzzy anti-bounded.
\end{theorem}
{\bf Proof.}$\;\;$
{\bf(i)$\;$} Since, $T$ is strongly fuzzy anti-continuous at
$x_0\,\in\;U$, for each $\epsilon\,>\,0\,$ there exists
$\,\delta\,>\,0\,$ such that   \[\nu_V\,(\;T(x)\,-\,T(x_0)\;,
\;\epsilon\;)\;\leq\;\nu_U\,(\;x\,-\,x_0\;,\;\delta\;)\]
Taking $\;y\,\in\;U\;$ and replacing $ \;x\;by\;x\,+\,x_0\,-\,y$, we get,\\
$\nu_V\,(\;T(x)\,-\,T(x_0)\;,\;\epsilon\;)\;\leq\;\nu_U\,(\;x\,
-\,x_0\;,\;\delta\;)\\ \Rightarrow\;\nu_V\,(\;T(x\,+\,x_0\,-\,y)\,
-\,T(x_0)\;,\;\epsilon\;)\;\leq\;\nu_U\, (\;x\,+\,x_0\,-\,y\,-\,x_0
\;,\;\delta\;)\\ \Rightarrow\;\nu_V\,(\;T(x)\,+\,T(x_0)\,-\,T(y)\,-
\,T(x_0)\;,\;\epsilon\;)\;\leq\;\nu_U\,(\,x\,-\,y\,,\,\delta\,)\\
\Rightarrow\;\nu_V\,(\;T(x)\,-\,T(y)\;,\;\epsilon\;)\;\leq\;
\nu_U\,(\,x\,-\,y\,,\,\delta\,)\\$
Since, $\,y\,\in\;U\,$ is arbitrary, $T$ is strongly fuzzy
anti-continuous on $U$.\\\\
{\bf(ii)$\;$}First we suppose that $T$ is strongly fuzzy anti-bounded.
Thus there exist a positive real number $M$ such that for all
$x\,\in\,U$ and for all $\epsilon\,\in\,R^+$,\\
\[\nu_V\,(\;T(x)\,,\,\epsilon\,)\;\leq\;\nu_U\,(\,x\,,
\,\frac{\epsilon}{M}\,)\hspace{7 cm}\]
\[i.e.,\;\nu_V\,(\;T(x)\,-\,T(\theta)\,,\,\epsilon\,)\;\leq\;
\nu_U\,(\,x\,-\theta\,,\,\frac{\epsilon}{M}\,)\hspace{6 cm}\]
\[i.e,\;\nu_V\,(\;T(x)\,-\,T(\theta)\,,\,\epsilon\,)\;\leq\;
\nu_U\,(\,x\,-\theta\,,\,\delta\,)\hspace{6 cm}\]
where  $\delta\;=\;\frac{\epsilon}{M}$.\\Thus $T$ is strongly
fuzzy anti-continuous at $\theta$  and hence $T$ is strongly
fuzzy anti-continuous on $U$.\\\\
Conversely, suppose that $T$ is strongly fuzzy anti-continuous
 on $U$. Using fuzzy anti-continuity of $T$ at $x\;=\;\theta$
 for $\epsilon\;=\;1$ there exist $\delta\;>\;0$ such that for
 all $x\;\in\;U$,\[\nu_V\,(\;T(x)\,-\,T(\theta)\;,\;1\;)\;\leq
 \;\nu_U\,(\,x\,-\,\theta\,,\,\delta\,)\] If $x\;\neq\;\theta\;
 \;$and $\;t\;>\;0\;.$ Putting $\;x\;=\;ut\;\\ \nu_V\,(\;T(x)\,,
 \,t\,)\;=\;\nu_V\,(\;uT(u)\,,\,t\,)\;=\;\nu_V\,(\;T(u)\,,\,1\,)\;
 \leq\;\nu_U\,(\;u\,,\,\delta\,)\;=
\;\nu_U\,(\;\frac{x}{t}\,,\,\delta\,)\;=\;\nu_U\,(\;x\,,\,
\frac{t}{M}\,)\;$, Where $\;M\;=\;\frac{1}{\delta}\,.\\$So,
$\;\nu_V\,(\;T(x)\,,\,t\;)\;\leq\;\nu_U\,(\;x\,,\,\frac{t}{M}\;)$\\
If $\;x\;\neq\;\theta\;$ and $\;t\;\leq\;0\;$ then $\;
\nu_V\,(\;T(x)\,,\,t\;)\;=\;1\;=\;\nu_U\,(\;x\,,\,\frac{t}{M}\;).\\$
If $\;x\;=\theta\;$ and $\;t\,\in\;R\,,\;$ then
$\;T(\theta_U)\;=\;\theta_V\;$ and
\[\nu_V\,(\;\theta_V\,,\,t\;)\;=\;\nu_U\,(\theta_U\,,
\,\frac{t}{M}\;)\;=\;0\;\;\;,\;if  \;\;\;t\,>\,0.\]
\[\nu_V\,(\;\theta_V\,,\,t\;)\;=\;\nu_U\,(\theta_U\,,\,
\frac{t}{M}\;)\;=\;1\;\;\;,\;if  \;\;\;t\,\leq\,1.\]
Hence $\;T\;$ is strongly fuzzy anti-bounded.

\begin{theorem}
Let $T\,:\,(\,U\;,\;A^\ast\,)\longrightarrow(\,V\;,\;B^\ast\,)$
be a linear operator. Then,\\
{\bf (i)$\;$} $T$ is weakly fuzzy anti-continuous on $U$ if $T$
is weakly fuzzy anti-continuous at a point $x_0\;\in\;U\,.$\\
{\bf (ii)$\;$} $T$ is weakly fuzzy anti-continuous if and only
 if $T$ is weakly fuzzy anti-bounded.
\end{theorem}
{\bf Proof.}
{\bf (i)$\;$} Since, $\;T\;$ is weakly fuzzy anti-continuous
at $\;x_0\;$ in $\;U\;$, for $\epsilon\,>\,0\;$ and
$\;\alpha\,\in\;(\,0\,,\,1\,)\;$ there  exist
$\;\delta\,=\,\delta\,(\alpha,\epsilon)\,>\,0\;$ such that
 $\forall\;x\,\in\,U$ \[\nu_U\,(\,x-x_0\,,\,\delta\,)\,\leq\,1-\alpha\;
 \;\Rightarrow\;\nu_V\,(\,T(x)\,-T(x_0)\,,\,\epsilon\,)\,
\leq\,1-\alpha.\;\]
Taking $\;y\,\in\,U\;$ and replacing $x$ by $\;x+x_0-y\;$ we get,\\
\[\nu_U\,(\,x+x_0-y-x_0\,,\,\delta\,)\,\leq\,1-\alpha\;\;
\Rightarrow\;\nu_V\,(\,T(x+x_0-y)\,-\,T(x_0)\,,\,
\epsilon\,)\,\leq\,1-\alpha\hspace{-1 cm}\]
\[i.e.,\;\nu_U\,(\,x-y\,,\,\delta\,)\,\leq\,1-\alpha\;\;\Rightarrow\;
\nu_V\,(\,T(x)\,+\,T(x_0)\,-\,T(y)\,-\,T(x_0)\,,\,
\epsilon\,)\,\leq\,1-\alpha\]
\[i.e.,\;\nu_U\,(\,x-y\,,\,\delta\,)\,\leq\,1-\alpha\;\;
\Rightarrow\;\nu_V\,(\,T(x)-T(y)\,,\,\epsilon\,)\,
\leq\,1-\alpha\hspace{5 cm}\]
Since, $\;y\,(\,\in\,U\,)\;$ is arbitrary it follows that $\;T\;$
is weakly fuzzy anti-continuous on $\;U\;.$ \\\\

{\bf (ii)$\;$} First we suppose that $\;T\;$ is fuzzy anti-bounded.
Thus for any $\;\alpha\,\in\;(\,0\,,\,1\,)\;$ there exist
$\;M_\alpha\,>\,0\;$ such that $\;\forall\,t\,\in\,R\;$ and
$\;\forall\,x\,\in\,U\;$ we have \[\nu_U(x,\frac{t}{M})\leq 1-
\alpha \Rightarrow \nu_V(T(x),t)\leq 1-\alpha\]
Therefore,
\[\nu_U(x-\theta,\frac{t}{M})\leq 1-\alpha \Rightarrow
\nu_V(T(x)-T(\theta),t)\leq 1-\alpha\hspace{3 cm}\]
\[i.e.,\,\nu_U(x-\theta,\frac{\epsilon}{M_\alpha})\leq 1-
\alpha \Rightarrow \nu_V(T(x)-T(\theta),\epsilon)\leq 1-
\alpha\hspace{5 cm}\]
\[i.e.,\nu_U(x-\theta,\delta)\leq 1-\alpha \Rightarrow
\nu_V(T(x)-T(\theta),\epsilon)\leq 1-\alpha\hspace{5 cm}\]
where $\delta=\frac{\epsilon}{M_\alpha}$\\
Thus, $T$ is weakly fuzzy anti-continuous at $x_0$ and
hence weakly fuzzy anti-continuous on $U$.\\
Conversely, suppose that $T$ is weakly fuzzy anti-continuous
 on $U$. Using continuity of $T$ at $\theta$ and taking
 $\epsilon=1$ we have for all $\alpha\in (0,1)$ there
 exists $\delta(\alpha,1)>0$ such that for all $x\in U,\\$
\[\nu_U(x-\theta,\delta)\leq 1-\alpha \Rightarrow
\nu_V(T(x)-T(\theta),1)\leq 1-\alpha\hspace{3.25 cm}\]
 \[i.e.,\,\nu_U(x,\delta)\leq 1-\alpha \Rightarrow
 \nu_V(T(x),1)\leq 1-\alpha.\hspace{6 cm}\]
If $x\neq \theta$ and $t>0$. Putting $x=\frac{u}{t}$
we have,\\
\[\nu_U(\frac{u}{t},\delta)\leq 1-\alpha \Rightarrow
\nu_V(T(\frac{u}{t}),1)\leq 1-\alpha \hspace{5 cm}\]
\[i.e.,\,\nu_U(u,t\delta)\leq 1-\alpha \Rightarrow
\nu_V(T(u),t)\leq 1-\alpha \hspace{6 cm}\]
\[i.e.,\,\nu_U(u,\frac{t}{M_\alpha})\leq 1-\alpha
\Rightarrow \nu_V(T(\frac{u}{t}),1)\leq 1-\alpha \hspace{7 cm}\]
 where $M_\alpha=\frac{1}{\delta(\alpha,1)}$
If $x\neq \theta$ and $t\leq 0$,  $\nu_U(x,
\frac{t}{M_\alpha})=\nu_V(T(x),t)=1$  for any $M_\alpha>0.\\$
If  $x=\theta$  then for $M_\alpha>0,\\$
\[\nu_U(x,\frac{t}{M_\alpha})=\nu_V(T(x),t)=0\;\;,\;if  \;\;\;t>0\]
\[\nu_U(x,\frac{t}{M_\alpha})=\nu_V(T(x),t)=1\;\;,\;if  \;\;\;t \leq 0\]
Hence,  $T$ is weakly fuzzy anti-bounded.

\begin{theorem}\label{t31}
Let $T\,:\,(\,U\;,\;A^\ast\,)\longrightarrow(\,V\;,\;B^\ast\,)$
 be a linear operator and $(\,U\;,\;A^\ast\,)$ and
  $(\,V\;,\;B^\ast\,)$ satisfies (vi) and (vii).
  Then $T$ is weakly fuzzy anti-bounded if and only if $T$ is fuzzy
anti-bounded with respect to $\alpha$-norms.
\end{theorem}
{\bf Proof.$\;\;$}
First we suppose that $T$ is weakly fuzzy anti-bounded.
Then for all $\alpha\,\in\;(0,1)$  there exist
$M_\alpha\,>\,0$ such that $\forall\,x\,\in\,U\,,\,t\,\in\,R\,$
we have\\
\[\nu_U(x,\frac{t}{M_\alpha})\,\leq\,1-\alpha\;\;\Rightarrow\;
\nu_V(T(x),t)\,\leq\,1-\alpha\;\]
Hence we get, $\nu_U(M_{\alpha}x,t)\,\leq\,1-\alpha\;\;
\Rightarrow\;\nu_V(T(x),t)\,\leq\,1-\alpha\;\\ i.e.,\;
\wedge\,\{\,\beta\,\in\,(0,1)\,:\,\left\| {\;M_\alpha\,x\;}
\right\|_{\,\beta}^{\,\ast}\;\leq\;t\,\}\;\leq\;1-\alpha \;
\Rightarrow\;\wedge\,\{\,\beta\,\in\,(0,1)\,:\,\left\| {\;T(x)\;}
\right\|_{\,\beta}^{\,\ast}\;\leq\;t\,\}\;\leq\;1-\alpha \;\\$
Now we show that\[\wedge\,\{\,\beta\,\in\,(0,1)\,:\,\left\|
{\;M_\alpha\,x\;} \right\|_{\,\beta}^{\,\ast}\;\leq\;t\,\}\;
\leq\;1-\alpha \; \Leftrightarrow\;\left\| {\;M_\alpha\,x\;}
\right\|_{\,\alpha}^{\,\ast}\;\leq\;t\]
If $x=\theta$ then the relation is obvious.\\
Suppose $x\neq \theta.\\$
Now, if
\begin{equation}\label{e51}
\wedge\,\{\,\beta\,\in\,(0,1)\,:\,\left\| {\;M_\alpha\,x\;}
\right\|_{\,\beta}^{\,\ast}\;\leq\;t\,\}\;<\;1-\alpha \;\;
then \;\; \left\| {\;M_\alpha\,x\;} \right\|_{\,\alpha}^{\,\ast}\;\leq\;t
\end{equation}
If $\;\wedge\,\{\,\beta\,\in\,(0,1)\,:\,\left\| {\;M_\alpha\,x\;}
\right\|_{\,\beta}^{\,\ast}\;\leq\;t\,\}\;=\;1-\alpha \;$
then there exists a decreasing sequence $\;\{\alpha_n\}_n\;$
 in $\,(0,1)\,$ such that $\;\alpha_n\;\longrightarrow\;\alpha\;$
 and $\;\left\| {\;M_\alpha\,x\;} \right\|_{\,\alpha}^{\,\ast}\;\leq\;t\;$
Then by Theorem \ref{th52} we have
\begin{equation} \label{e52}
\;\left\| {\;M_\alpha\,x\;} \right\|_{\,\alpha}^{\,\ast}\;\leq\;t
\end{equation}
From (\ref{e51}) and (\ref{e52}) we get
\begin{equation}\label{e53}
\wedge\,\{\,\beta\,\in\,(0,1)\,:\,\left\| {\;M_\alpha\,x\;}
 \right\|_{\,\beta}^{\,\ast}\;\leq\;t\}\;\leq\;1-\alpha \;
 \Rightarrow\;\left\| {\;M_\alpha\,x\;} \right\|_{\,\alpha}^{\,\ast}\;\leq\;t
\end{equation}
Next we suppose that $\;\left\| {\;M_\alpha\,x\;}
 \right\|_{\,\alpha}^{\,\ast}\;\leq\;t.\\$
If $\;\left\| {\;M_\alpha\,x\;} \right\|_{\,\alpha}^{\,\ast}\;<\;t\;$
then $\;\nu_U(\,M_\alpha\,x\,,\,t\,)\;\leq\;1-\alpha.\;$ i.e.,
\begin{equation}\label{e54}
 \;\wedge\,\{\,\beta\,\in\,(0,1)\,:\,\left\| {\;M_\alpha\,x\;}
 \right\|_{\,\beta}^{\,\ast}\;\leq\;t\}\;\leq\;1-\alpha
\end{equation}
If $\;\left\| {\;M_\alpha\,x\;} \right\|_{\,\alpha}^{\,\ast}\;=\;t\;$
i.e.,$\;\wedge\,\{\,s\,:\,\nu_U(M_\alpha\,x\,,\,s\,)\;
\leq\;1-\alpha\}\;=\;t\;$ then there exist an increasing
sequence $\{s_n\}_n$ in $\mathbf{R}^+$ such that
$\;s_n\,\longrightarrow\,t\,$ and
\[\nu_U(\,M_\alpha\,x\,,\,s_n\,)\;\leq\;1-\alpha\;\Rightarrow\;
\mathop{\lim }\limits_{n\, \to\,\;\infty}\, \nu_U(\,M_\alpha\,x\,,
\, s_n\,)\;\leq\;1-\alpha\]
\[\Rightarrow\;\nu_U(\,M_\alpha\,x\,,\,\mathop{\lim }\limits_{n\, \to\,
\;\infty}\,s_n\,)\;\leq\;1-\alpha\]
\[\Rightarrow\;\nu_U(\,M_\alpha\,x\,,\,t\,)\;\leq\;1-\alpha\hspace{1.25 cm}\]
\[\Rightarrow\;\wedge\,\{\,\beta\,\in\,(0,1)\,:\,
\left\| {\;M_\alpha\,x\;} \right\|_{\,\beta}^{\,\ast}\;
\leq\;t\}\;\leq\;1-\alpha\hspace{-2.5 cm}\]
Hence from (\ref{e54}) it follows that
\begin{equation}\label{e55}
\;\left\| {\;M_\alpha\,x\;} \right\|_{\,\alpha}^{\,\ast}\;\leq\;t\;
\Rightarrow\; \wedge\,\{\,\beta\,\in\,(0,1)\,:\,\left\| {\;M_\alpha\,x\;} \right\|_{\,\beta}^{\,\ast}\;\leq\;t\}\;\leq\;1-\alpha
\end{equation}
From (\ref{e53}) and (\ref{e55}) we have
\begin{equation}\label{e56}
\wedge\,\{\,\beta\,\in\,(0,1)\,:\,\left\| {\;M_\alpha\,x\;}
 \right\|_{\,\beta}^{\,\ast}\;\leq\;t\}\;\leq\;1-\alpha\;
 \Leftrightarrow\;\left\| {\;M_\alpha\,x\;} \right\|_{\,\alpha}^{\,\ast}\;\leq\;t
\end{equation}
In the similar way we can show that
\begin{equation}\label{e57}
\wedge\,\{\,\beta\,\in\,(0,1)\,:\,\left\| {\;T(x)\;}
\right\|_{\,\beta}^{\,\ast}\;\leq\;t\}\;\leq\;1-\alpha\;
\Leftrightarrow\;\left\| {\;T(x)\;} \right\|_{\,\alpha}^{\,\ast}\;\leq\;t
\end{equation}
From (\ref{e56}) and (\ref{e57}) we have $\;\;
\nu_U(\,M_\alpha\,x\,,\,t\,)\;\leq\;1-\alpha\;\Rightarrow\;
\nu_V(\,T(x)\,,\,t\,)\;\leq\;1-\alpha\\$
Then \[\left\| {\;M_\alpha\,x\;} \right\|_{\,\alpha}^{\,\ast}\;
\leq\;t\;\Rightarrow\;\left\| {\;T(x)\;} \right\|_{\,\alpha}^{\,\ast}\;\leq\;t\]
This implies that \[\left\| {\;T(x)\;} \right\|_{\,\alpha}^{\,\ast}\;
\geq\;\left\| {\;M_\alpha\,x\;} \right\|_{\,\alpha}^{\,\ast}\]
Conversely, suppose that $\;\forall\;\alpha\,\in\,(0,1),\;
\exists\;M_\alpha\,>\,0\;$ such that $\;\forall\;x\,In\,U,
$\[\left\| {\;T(x)\;} \right\|_{\,\alpha}^{\,\ast}\;\geq\;
\left\| {\;M_\alpha\,x\;} \right\|_{\,\alpha}^{\,\ast}\]
Then for $x\,\neq\,\theta\;\;and\;\;\forall\;t\,>\,0,$
\[\left\| {\;M_\alpha\,x\;} \right\|_{\,\alpha}^{\,\ast}\;
\leq\;t\;\Rightarrow\;\left\| {\;T(x)\;} \right\|_{\,\alpha}^{\,\ast}\;\leq\;t\]
i.e.,\[\wedge\,\{\,s\;:\;\nu_U(\,M_\alpha\,x\,,\,s\,)\;\leq\;
1-\alpha\;\}\;\leq \;t\;\Rightarrow\; \wedge\,\{\,s\;:\;
\nu_V(\,T(x)\,,\,s\,)\; \leq\;1-\alpha\;\}\;\leq \;t\]
In the similar way as above we can show that
\[\wedge\,\{\,s\;:\;\nu_U(\,M_\alpha\,x\,,\,s\,)\;\leq\;1-
\alpha\;\}\;\leq \;t\;\Leftrightarrow\;\nu_U(\,M_\alpha\,x\,,
\,t\,)\;\leq\;1-\alpha\]and
\[\wedge\,\{\,s\;:\;\nu_U(\,T(x)\,,\,s\,)\;\leq\;1-\alpha\;\}\;
\leq \;t\;\Leftrightarrow\;\nu_U(\,T(x)\,x\,,\,t\,)\;\leq\;1-\alpha\]
Thus we have \[\nu_U(\,x\,,\,\frac{t}{M_\alpha})\;\leq\;1-\alpha\;
\Rightarrow\;\nu_V(\,T(x)\,,\,t\,)\;\leq\;1-\alpha\;,\;\forall\,x\,\in\,U\]
If $\;x\,\neq\,\theta,\;t\,\leq\,0\;$ and if $\;x\,=\,
\theta,\;t\,>\,0\;$ then the above relation is obvious.
Hence the proof.

\begin{theorem}
Let, $T\,:\,(\,U\;,\;A^\ast\,)\longrightarrow(\,V\;,\;B^\ast\,)$
 be a linear operator and $(\,U\;,\;A^\ast\,)$ and $(\,V\;,\;B^\ast\,)$
 satisfies (vi) and (vii). If $U$ is finite dimensional then $T$
 is weakly fuzzy anti-bounded.
\end{theorem}
{\bf Proof.$\;\;$}
Since, $(\,U\;,\;A^\ast\,)$ and $(\,V\;,\;B^\ast\,)$ satisfies (vi),
 we may suppose that $\{\;\left\| {\;\cdot\;}
 \right\|_{\,\alpha }^{\,\ast}\;:\;\alpha\,\in\,(\,0\,,\,1\,)\}$
 is ascending family of fuzzy $\alpha$-anti-norms.\\\\Since $T$
 is of finite dimension, $T\,:\,(\,U\;,\;A^\ast\,)
 \longrightarrow(\,V\;,\;B^\ast\,)\;$ is bounded linear operator
 for each $\alpha\,\in\,(\,0\,,\,1\,)\,.\;$ Thus by Theorem
 \ref{t31} it follows that $T$ is weakly fuzzy anti-bounded.

\section{Future Work}
 In our next paper we shall try to develop the concept of fuzzy anti-bounded linear functionals and their properties.


\begin{thebibliography}{0}
\bibitem{Samanta1} Iqbal H. Jebril  and T.K. Samanta  ,
\textit{Fuzzy anti-normed linear space}, Journal of
mathematics and Technology, February, 2010.
\bibitem{Bag1}T. Bag and S.K. Samanta ,
\textit{ Finite Dimensional Fuzzy Normed Linear Spaces}, The Journal
of Fuzzy Mathematics Vol. 11 $(\,2003\,)$ 687 - 705.
\bibitem{Bag2} T. Bag and S.K. Samanta ,
\textit{Fuzzy bounded linear operators} , Fuzzy Sets and Systems 151
$(\,2005\,)$ 513 - 547.
\bibitem{Shih-chuan} S.C. Cheng  and J.N. Mordeson  ,
\textit{Fuzzy Linear Operators and Fuzzy Normed Linear
Spaces} , Bull. Cal. Math. Soc.86 $(\,1994\,)$ 429 - 436.
\bibitem{dinda1} Bivas Dinda, T.K. Samanta and Iqubal H. Jebril ,
\textit{Fuzzy Anti-norm and Fuzzy $\alpha$-anti-convergence} ,
(Comunicated)
\bibitem{dinda} Bivas Dinda  and T.K. Samanta  ,
\textit{Intuitionistic Fuzzy Continuity and Uniform Convergence} ,
Int. J. Open Problems Compt.Math., Vol 3, No. 1 $(\,2010\,)$ 8 - 26.
\bibitem{Felbin1}C. Felbin  ,
\textit{The completion of fuzzy normed linear space}, Journal
 of mathmatical analysis and application 174(2) $(\,1993\,)$ 428-440.
\bibitem{Samanta}T.K. Samanta  and Iqbal H. Jebril ,
\textit{Finite dimentional intuitionistic fuzzy normed linear space},
Int. J. Open Problems Compt. Math., Vol 2, No. 4 $(\,2009\,)$ 574-591.
\bibitem{Katsaras}A.K. Katsaras  ,
\textit{ Fuzzy topological vector space}, Fuzzy Sets and
Systems 12 $(\,1984\,)$ 143 - 154.
\bibitem{Schweizer}B. Schweizer , A. Sklar,
\textit{Statistical metric space}, Pacific journal of
 mathhematics 10 $(\,1960\,)$ 314-334.
\bibitem{Vijayabalaji}S. Vijayabalaji, N. Thillaigovindan, Y.B. Jun
\textit{Intuitionistic Fuzzy n-normed linear space} ,
Bull. Korean Math. Soc. 44 $(\,2007\,)$ 291 - 308.
\bibitem{Kramosil}O. Kramosil, J. Michalek   ,
\textit{ Fuzzy metric and statisticalmetric spaces},
 Kybernetica 11 $(\,1975\,)$ 326 - 334.
\bibitem{zadeh} L.A. Zadeh
\textit{Fuzzy sets}, Information and control 8 $(\,1965\,)$ 338-353.
\end{thebibliography}
\end{document}